\documentclass{article}
\usepackage{arxiv}

\usepackage[backend=biber,
sorting=none,
citestyle=numeric-comp]{biblatex}
\usepackage[utf8]{inputenc} 
\usepackage[T1]{fontenc}    
\usepackage{hyperref}       
\usepackage{url}            
\usepackage{booktabs}       
\usepackage{amsfonts}       
\usepackage{nicefrac}       
\usepackage{microtype}      
\usepackage{lipsum}		
\usepackage{graphicx}
\usepackage{doi}
\usepackage{
            amssymb,
            url,
            float,
            caption,
            subcaption,
            cleveref,
            subcaption,
            enumitem
            }

\usepackage{listings}
\usepackage{xcolor}

\usepackage{hyperref}
\hypersetup{
    colorlinks=true,
    linkcolor=blue,
    filecolor=magenta,      
    urlcolor=cyan,
    pdftitle={Overleaf Example},
    pdfpagemode=FullScreen,
    }

\definecolor{codegreen}{rgb}{0,0.6,0}
\definecolor{codegray}{rgb}{0.5,0.5,0.5}
\definecolor{codepurple}{rgb}{0.58,0,0.82}
\definecolor{backcolour}{rgb}{0.95,0.95,0.92}

\lstdefinestyle{mystyle}{
    backgroundcolor=\color{backcolour},   
    commentstyle=\color{codegreen},
    keywordstyle=\color{magenta},
    numberstyle=\tiny\color{codegray},
    stringstyle=\color{codepurple},
    basicstyle=\ttfamily\footnotesize,
    breakatwhitespace=false,         
    breaklines=true,                 
    captionpos=b,                    
    keepspaces=true,                 
    numbers=left,                    
    numbersep=5pt,                  
    showspaces=false,                
    showstringspaces=false,
    showtabs=false,                  
    tabsize=2
}

\title{A Python Benchmark Functions Framework for Numerical Optimisation Problems}
\author{ \href{https://orcid.org/0000-0002-2883-706X}{\includegraphics[scale=0.06]{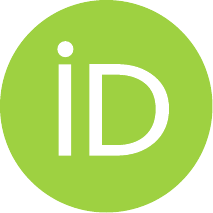}\hspace{1mm}Luca Baronti}\thanks{School of Computer Science, University of Birmingham, Birmingham, UK}, \href{https://orcid.org/0000-0002-5623-7491}{\includegraphics[scale=0.06]{orcid.pdf}\hspace{1mm}Marco Castellani}\thanks{Department of Mechanical Engineering, University of Birmingham, Birmingham, UK}}

\begin{document}
\maketitle

\begin{abstract}
This work proposes a framework of benchmark functions designed to facilitate the creation of test cases for numerical optimisation techniques. 
The framework, written in Python 3, is designed to be easy to install, use, and expand. 
The collection includes some of the most used multi-modal continuous functions present in literature, which can be instantiated using an arbitrary number of dimensions.
Meta-information of each benchmark function, like search boundaries and position of known optima, are included and made easily accessible through class methods.
Built-in interactive visualisation capabilities, baseline techniques, and rigorous testing protocols complement the features of the framework. The framework can be found here: \url{https://gitlab.com/luca.baronti/python_benchmark_functions}

\end{abstract}











\section{Motivation}\label{sec:motivations}
Numerical optimisation problems are often used to model real-world situations \cite{tsai2014optimization} where, given a certain function, we are interested in finding one or more minima or maxima.
Different optimisation techniques are usually evaluated in terms of computation time or number of sampling in the search space. Unfortunately, an optimal technique that outperforms another in every problem can not \cite{wolpert1997no} exists. Therefore, identifying an appropriate subset of problems (i.e. functions) where a novel technique works particularly well is of primary importance.

Despite some functions having become standard benchmark cases \cite{pham2015comparative} in the field, they are often implemented from scratch by the practitioners, involving unnecessary extra work and the risk of introducing implementation errors \cite{posypkin2017implementation} which may be hard to identify. Minimal differences in the actual implementation may result in numerical differences in the function evaluation. This can lead to erroneous comparison with alternative techniques when early stopping criteria, based on mismatching information of local optima, are used.

Collections of benchmark functions \cite{averick1991minpack,ali2005numerical,pham2014benchmarking,vanaret2020certified,pohlheim2007examples,jamil2013literature} have been published, but unfortunately those efforts haven't been coupled with off-the-shelf implementations of the described functions. Attempts to bridge the gap between the formal description of benchmark functions and usable frameworks have been proposed in the literature. MVF \cite{adorio2005mvf} is a C library\footnote{A version of the C library MVF can be found here: \url{https://gitlab.com/luca.baronti/mvf}} that provides the implementation of a number of different benchmark functions. Tests suites for optimisation and linear algebra solvers have also been proposed, like CUTEst \cite{gould2015cutest}. More recently, a C++ library has been proposed \cite{posypkin2017implementation} which includes a large number of functions along with some meta-data information.

The proposed framework\footnote{The framework can be found here: \url{https://gitlab.com/luca.baronti/python_benchmark_functions}} differs from many solutions already available for:
\begin{enumerate}[label=(\roman*)]
    \item its flexibility, simplicity of use and expansion: the framework can be easily installed with no compilation steps, using the PyPI system. Most functions can be instantiated in an arbitrary number of dimensions, using them requires only the minimal number of operations necessary (see \cref{sec:examples}) and new functions can be easily incorporated (see \cref{sec:arch});
    \item including suggested search boundaries, the position of the known global minimum/maximum along relevant local minima/maxima, useful for testing multi-optima techniques and making them readily available as class methods;
    \item including advanced visualisation capabilities: functions can be easily plotted as lines, surfaces or heatmaps. Optionally, a set of points can be also visualised along the function, allowing the practitioner to follow the progress of their technique in real time (see \cref{sec:visualisation});
    \item including useful meta-information: the description of a function, the BibTex reference to the original paper as well as its definition (in \LaTeX{}) are examples of information that can be directly accessed through class methods;
\end{enumerate} 
The inclusion of two baseline search methods and a local optimum tester complement the set of features offered by the framework.






\section{Framework Description}
\label{sec:framework_description}
\begin{figure}[ht]
\begin{subfigure}[b]{.3\linewidth}
\includegraphics[width=\linewidth]{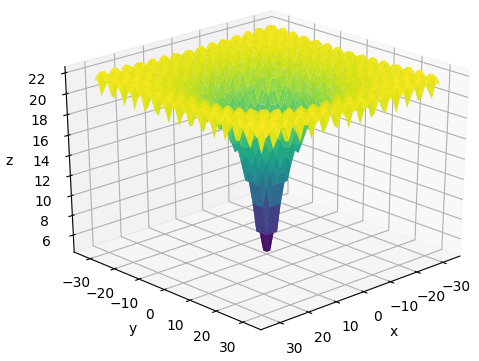}
\caption{Ackley\cite{ackley2012connectionist}}
\end{subfigure}%
\begin{subfigure}[b]{.3\linewidth}
\includegraphics[width=\linewidth]{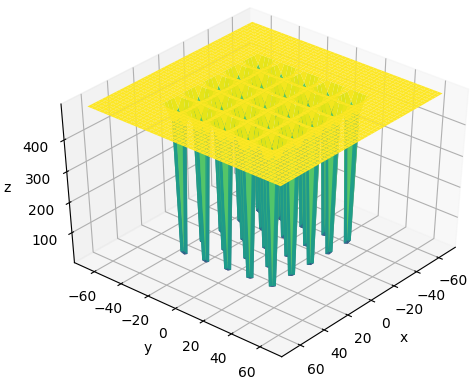}
\caption{De Jong 5\cite{molga2005test}}
\end{subfigure}%
\begin{subfigure}[b]{.3\linewidth}
\includegraphics[width=\linewidth]{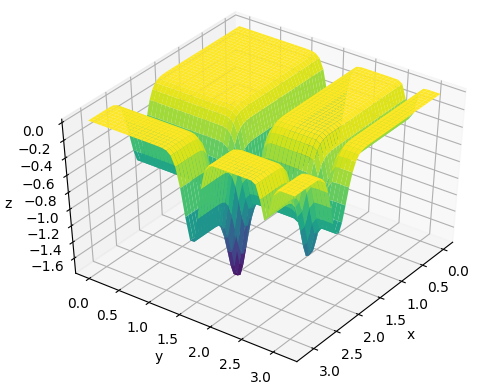}
\caption{Michalewicz\cite{vanaret2020certified}}
\end{subfigure}

\begin{subfigure}[b]{.3\linewidth}
\includegraphics[width=\linewidth]{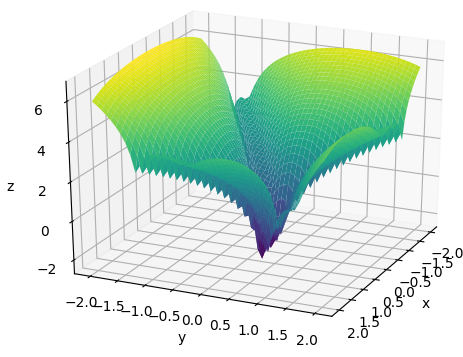}
\caption{Picheny, Goldstein and Price\cite{picheny2013benchmark}}
\end{subfigure}%
\begin{subfigure}[b]{.3\linewidth}
\includegraphics[width=\linewidth]{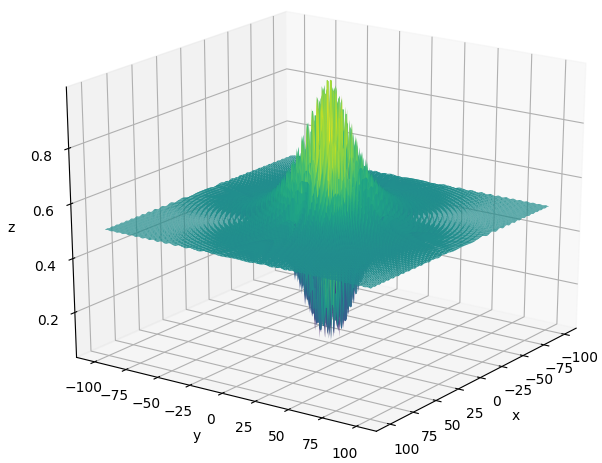}
\caption{Schaffer N2\cite{mishra2006some}}
\end{subfigure}%
\begin{subfigure}[b]{.3\linewidth}
\includegraphics[width=\linewidth]{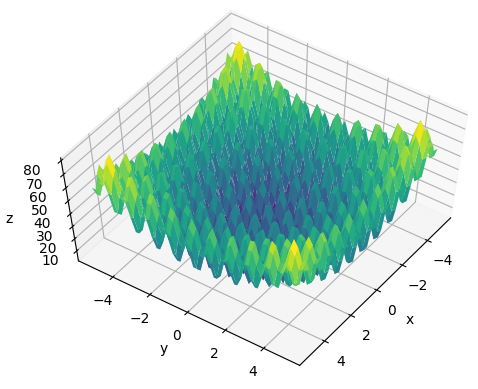}
\caption{Rastring\cite{pohlheim2007examples}}
\end{subfigure}%

\begin{subfigure}[b]{.3\linewidth}
\includegraphics[width=\linewidth]{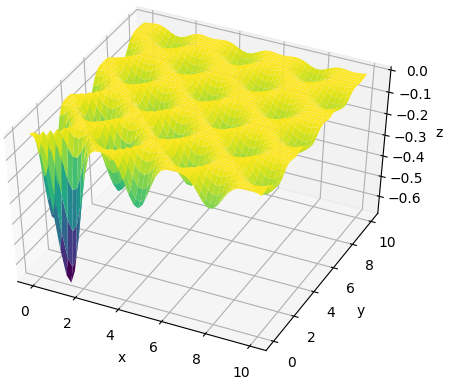}
\caption{Keane\cite{vanaret2020certified}}
\end{subfigure}%
\begin{subfigure}[b]{.3\linewidth}
\includegraphics[width=\linewidth]{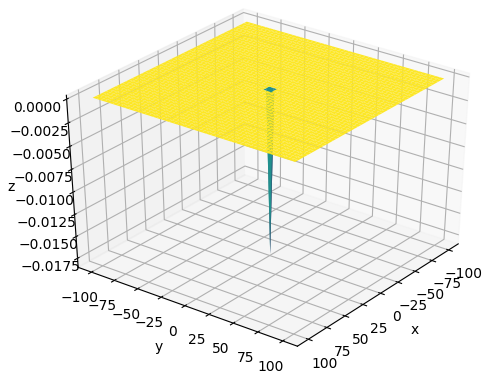}
\caption{Easom\cite{chung1998caep}}
\end{subfigure}%
\begin{subfigure}[b]{.3\linewidth}
\includegraphics[width=\linewidth]{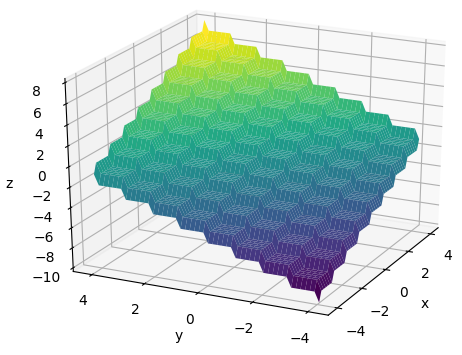}
\caption{De Jong 3\cite{de1975analysis}}
\end{subfigure}%

\caption{Example of benchmark functions included in the library.}\label{fig:functions}
\end{figure}

The proposed framework is a collection of multi-modal benchmark functions designed to be used for testing numerical optimisation techniques on multi-dimensional continuous problems. This initial version of the library includes 20 of the most common functions (a sample is visible in \Cref{fig:functions}). Most of them can be instantiated using an arbitrary number of dimensions. The library has been designed to be easy to use (see \Cref{sec:examples}) and easily expandable (see \Cref{sec:arch}). Every function is coupled with metadata (see \Cref{sec:funcs}) that are accessible out-of-the-box and generally necessary to the practitioners. Finally, every function can be plotted in an interactive view using just one line of code (see \Cref{sec:visualisation}).

\subsection{Framework Architecture}
\label{sec:arch}
The library is divided in a \emph{benchmark\_functions.py} and a \emph{functions\_info\_loaders.py} files. The former is the main file that contains the abstract class \emph{BenchmarkFunction} that every benchmark function extends, as well as the implementation of each function available. The latter is an auxiliary file, which contains some simple data structures and the routines used to load the functions metadata. Every function metadata is present in a JSON file in the \emph{functions\_info} directory. Generally, the optima of a function differs w.r.t. the number of dimensions. For this reason, the JSON schema groups the optima according to their dimension, using a special character (\emph{*}) for the optima that can be analytically proven to be dimensional invariant. The JSON schema has been designed to also support parametric functions.

The metadata format is completely transparent to the end user, which accesses the information through class methods. This architecture is designed to make interventions on specific functions metadata (e.g. the modification of the known optima) as local as possible, and the expansion of the library with new functions as easy as possible. The library can be easily installed with the PyPI system:
\begin{lstlisting}[language=bash]
$ pip3 install benchmark_functions
\end{lstlisting}
A novel benchmark function can be implemented as a class that extends \emph{BenchmarkFunction} overriding the \emph{\_evaluate()} method, which returns the function's value at a given point, and adding the metadata in a dedicated JSON file in the proper directory.

A Continuous Integration (CI) testing framework has been implemented to mitigate the inclusion of bugs in future updates. This framework performs some sanity checks, including a test on the local optima provided by the library. This test samples a large number of points within a very small hyper-sphere centred in the putative optimum, verifying that no better solutions, within a given threshold $\epsilon$, have been found this way. All the optima for the functions provided have been tested for $\epsilon=10^{-6}$. Following the CI paradigm, this test routine is automatically performed every time a new version is pushed to the repository.

\subsection{Framework Functionalities}
\label{sec:funcs}
To use a function it is sufficient to create an instance of it and call the instance on a list (see \cref{sec:examples}). The library includes functions defined for a fixed or an arbitrary number of dimensions $N$. In the latter case, $N$ can be specified in the constructor.  In case of parametric functions (e.g. Ackley \cite{ackley2012connectionist}) those parameters can also be set in the constructor. Maximisation and minimisation problems are fundamentally equivalent. Traditionally, benchmark functions have been designed as minimisation problems, whilst many optimisation algorithms use a \emph{score} function that needs to be maximised. To facilitate the use in these cases, a \emph{opposite} flag can be set to \emph{True} to instantiate the opposite version of the function (see \cref{sec:optim_example}). 

Basic function information can be accessed with relevant methods (e.g. \emph{name()}, \emph{n\_dimensions()}) whilst the method \emph{suggested\_bonds()} returns the boundaries where the function is usually evaluated in literature. The \emph{minimum()} and \emph{maximum()} methods provide the best known approximation of optimum of the function. Some benchmark functions (e.g. Ackley) are designed to have numerous deceptive local optima with small attraction basins, others (e.g. Michalewicz) have few distinct strict local optima. In the latter case, information about those local optima are as important as the information about the global optimum for multi-optima optimisation techniques. In contrast with similar libraries, this framework also provides this information with relevant methods (e.g. \emph{n\_minima()}, \emph{minima()}, \emph{saddle\_points()}).

A \LaTeX{} definition of the function can be accessed with the \emph{definition()} method, whilst the \emph{reference()} method provides a BibTex entry for the paper that first introduced it\footnote{In some cases it has not been possible to trace the exact origin of certain functions. In these cases, the reference is either empty or it is the oldest known source that mentioned the function.}. These methods are meant to facilitate the inclusion of a function in a publication. 

The \emph{show()} method displays an interactive plot (see \cref{sec:visualisation}) of the function (provided $N\in\{1,2\}$). For $N=2$ an heatmap representation can alternatively be chosen, setting \emph{asHeatMap} to \emph{True} (see \cref{fig:heatmapWpoints}). Optionally, a set of points can be passed to the method to be plotted along the function. Alternative boundaries can be chosen with the \emph{bounds} parameter. This allows us to visualise a specific sub-region of the function.

Two baseline algorithms are also included: \emph{minimum\_grid\_search()} and \emph{minimum\_random\_search()}. The latter performs a simple random sampling of \emph{n\_samples} points whilst the former performs a grid search of $($\emph{n\_edge\_points}$+1)^N$ points. In both cases the sampling is performed within the suggested boundaries (by default) and the sampled point of minimum value is returned. These are not meant to be used as proper optimisation algorithms, rather they are included just to provide out-of-the-box baseline techniques useful to compare the performances of a novel optimisation algorithm.

\section{Illustrative Examples}
\label{sec:examples}
To use a function from the collection (e.g. Schwefel \cite{schwefel1981numerical}) it is sufficient to instantiate its class from the library:
\begin{lstlisting}[language=Python]
>>> import benchmark_functions as bf
>>> 
>>> func = bf.Schwefel(n_dimensions=4)
>>> point = [25, -34.6, -112.231, 242]
>>> func(point) # function's value at the given point
-129.38197657025287
\end{lstlisting}
Most functions implemented can be instantiated with an arbitrary number of dimensions. This can be set with the \emph{n\_dimensions} optional parameter. If the number of dimensions is not specified a default value (generally $N=2$) will be used. Parameters required by parametric functions can be set in the constructor, otherwise default values will be taken in these cases. Some functions are only defined for 2 dimensions (e.g. Easom) in these cases no \emph{n\_dimensions} parameter is accepted.
\subsection{Visualisation}\label{sec:visualisation}
If the number of dimensions is either 1 or 2 it is possible to plot the benchmark function in an interactive widget.
The resulting plot is either a 3D surface (when $N=2$) or a simple 2D graph plot ($N=1$). 
If the function is defined in 2 dimensions, it is also possible to plot it as an heatmap setting the function parameter \emph{asHeatMap=True} as follows:
\begin{lstlisting}[language=Python]
>>> func = bf.Schwefel(n_dimensions=2)
>>> func.show(asHeatMap=True)
\end{lstlisting}
By default, the function will be shown within its suggested boundaries. It is possible to pass custom boundaries, for visualisation purposes, using the parameter \emph{bounds}.

\begin{figure}[ht]
    \centering
    \begin{subfigure}[b]{.49\textwidth}
    \includegraphics[width=\linewidth]{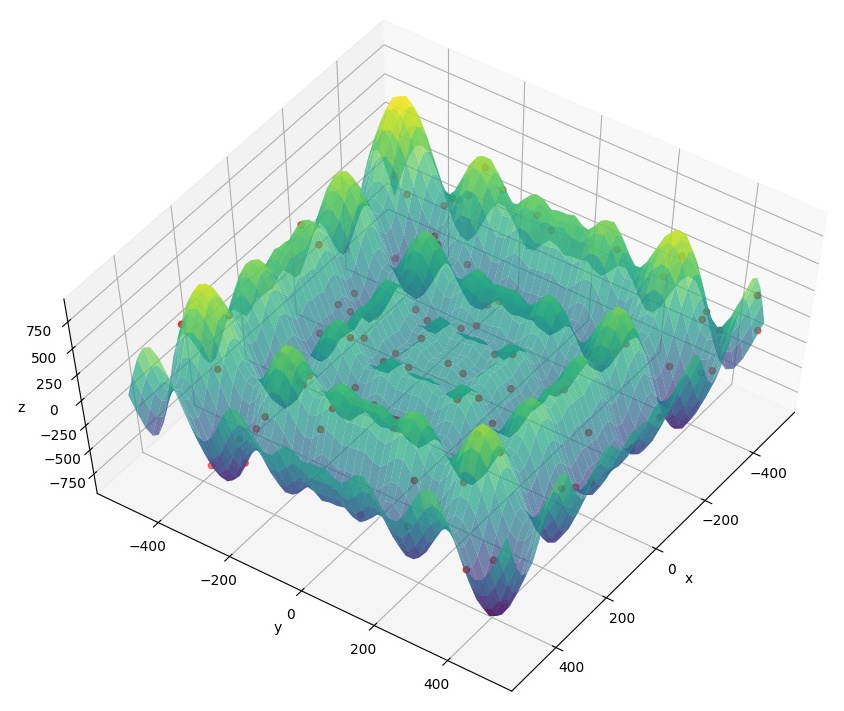}
    \caption{3D Plot with points}
    \label{fig:plotWpoints}
    \end{subfigure}
    \begin{subfigure}[b]{.49\textwidth}
    \includegraphics[width=\linewidth]{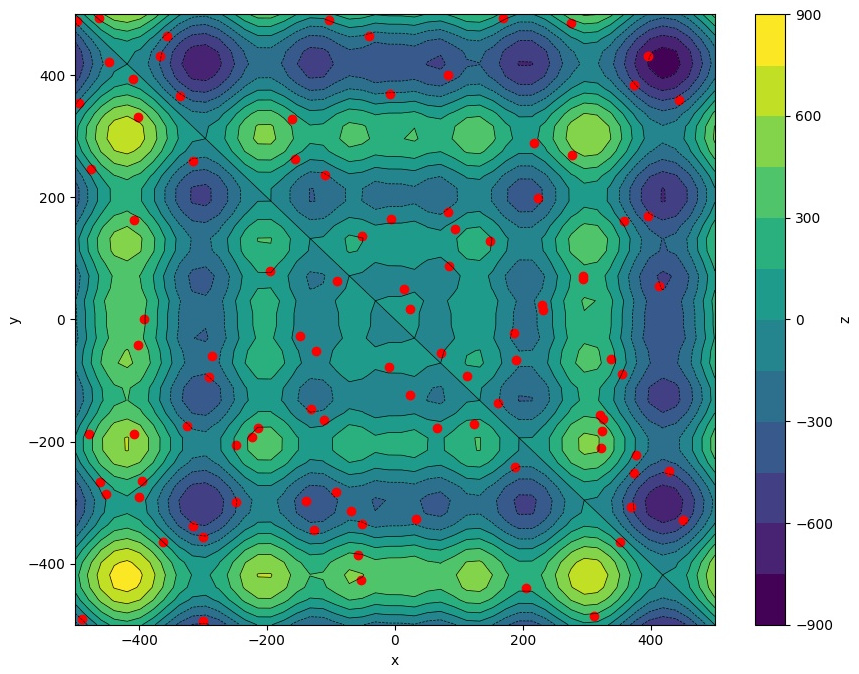}
    \caption{Heatmap with points}
    \label{fig:heatmapWpoints}
    \end{subfigure}
    \caption{Example of visualisation of the Schwefel function with a set of random points, uniformly sampled within the suggested boundaries.}
    \label{fig:exampleWpoints}
\end{figure}

A set of points can be visualised over the function's surface (or on the heatmap) passing a list of (1D or 2D) points to the parameter \emph{showPoints}. For instance, following the previous snippet, 100 random points can be generated and displayed as follows:
\begin{lstlisting}[language=Python]
>>> import numpy as np
>>> 
>>> lb, ub = func.suggested_bounds()
>>> random_points = [np.random.uniform((lb, lb), (ub, ub), 2) for _ in range(100)]
>>> func.show(showPoints=random_points)
\end{lstlisting}
This generates the picture shown in \Cref{fig:plotWpoints}. The points can be displayed on the heatmap in similar fashion:
\begin{lstlisting}[language=Python]
>>> func.show(showPoints=random_points, asHeatMap=True)
\end{lstlisting}
producing the picture visible in \Cref{fig:heatmapWpoints}.
\subsection{Example of Use with an Optimisation Algorithm}
\label{sec:optim_example}
Here is provided an example of out-of-the-box use of the library in the context of testing a  numerical optimisation technique.
In this example we want to test the Bees Algorithm \cite{baronti2020analysis} against the De Jong 5 function, in order to find the point of minimum value. The first step is importing the library and the optimisation algorithm\footnote{A Python version of the Bees Algorithm is available here: \url{https://gitlab.com/bees-algorithm/bees_algorithm_python}}:
\begin{lstlisting}[language=Python]
>>> import benchmark_functions as bf
>>> from bees_algorithm import BeesAlgorithm
\end{lstlisting}
Then the function is instantiated and the suggested boundaries are used to restrict the scope of the search. The implementation of the Bees Algorithm used works on maximisation problems, whilst we are interested in finding the minimum of the function. This can be solved just by passing the \emph{opposite=True} flag to the function's constructor.
\begin{lstlisting}[language=Python]
>>> func = bf.DeJong5(opposite=True)
>>> lb, ub = func.suggested_bounds()
>>> alg = BeesAlgorithm(score_function=lambda x: func(x), range_min=lb, range_max=ub)
\end{lstlisting}
Here the optimisation is performed:
\begin{lstlisting}[language=Python]
>>> alg.performFullOptimisation(max_iteration=3000)
\end{lstlisting}
Once the optimisation is ended, we can recover the best solution found by the optimiser:
\begin{lstlisting}[language=Python]
>>> local_minimum = alg.best_solution.values
>>> print(local_minimum)
[-31.978333633653907, -31.978334898791008]
\end{lstlisting}
One of the advantages of using this library is that this solution 
\begin{lstlisting}[language=Python]
>>> print((func(local_minimum), local_minimum))
(0.9980038377944496, [-31.978333633653907, -31.978334898791008])
\end{lstlisting}
can  be directly compared with the best approximation available in literature:
\begin{lstlisting}[language=Python]
>>> func = bf.DeJong5()
>>> print(func.minimum())
(0.9980038377944496, [-31.978333625355454, -31.978335021953196])
\end{lstlisting}
As it is possible to see, although the solution found by the algorithm is slightly different from the solution present in the library, in both cases the function's value is equivalent.

\section{Impact of the Proposed Framework}
\label{sec:impact}
The main goal of this library is to provide a powerful and flexible benchmark framework for testing optimisation techniques. The framework gives the practitioners easy access to an ever growing number of curated benchmark functions, removing the time and risks involved in the development of custom solutions. Moreover, the standardisation of the test framework directly promotes the consistency and accuracy in the results among different publications. The direct availability of search boundaries, optima positions and their values further helps the testing of a novel technique against a large number of functions, whilst the easy access to functions meta-data streamline their inclusion in a publication.

The in-built visualisation capabilities (see \cref{sec:visualisation}) provide an easy and versatile aid in the design of novel optimisation techniques. Intermediate solutions of iterative optimisation techniques can be directly visualised in the interactive plot, providing useful insights \cite{baronti2020analysis} on the behaviour of the algorithm. This is intended to help the practitioner in early identifying potential problems with a novel technique, but it is also an useful tool for presenting the algorithm in live demos as well as for didactic purposes. 

Finally, the library's architecture has been designed with a collaborative approach in mind. New interesting functions and new optima can be easily added with local interventions, and new additions are monitored by the CI test suite provided. 
\section{Future Works}
\label{sec:future}

The collection of benchmark functions included is already varied enough to cover many test cases. However, expanding the collection with other functions used in literature is one of the short term goals for the future developments. The inclusion of the gradient and Hessian definitions for each function is also a priority. This will expand the scope of the framework, making it available also for gradient-based optimisation techniques. Finally, adding a structure to the collection (e.g. grouping the functions by class, or establishing a taxonomy) can guide the practitioner in choosing a subset of functions fitting the needs of their specific test case.

\section*{Declaration of competing interest}
The author wish to confirm that there are no known conflicts of interest associated with this publication and there has been no significant financial support for this work that could have influenced its outcome.

\printbibliography




\end{document}